\documentclass{gtart_h}  

\def\ifplaintex{\expandafter\ifx\csname documentclass\endcsname\relax}

\ifplaintex 
\else
\fi

\expandafter\ifx\csname beginpicture\endcsname\relax
\expandafter\ifx\csname documentclass\endcsname\relax
\input pictex \else
\input prepictex \input pictex \input postpictex \fi\fi

\def\gt{{\mathsurround=0pt\it $\cal G\mskip-2mu$eometry \&\ 
$\cal T\!\!$opology}}        

\def\gtp{{\mathsurround=0pt\it $\cal G\mskip-2mu$eometry \&\ 
$\cal T\!\!$opology $\cal P\!$ublications}}  

\def\lognumber#1{\def\thelognumber{#1}}
\def\volumenumber#1{\def\thevolumenumber{#1}}
\def\papernumber#1{\def\thepapernumber{#1}}
\def\volumeyear#1{\def\thevolumeyear{#1}}

\def\pagenumbers#1#2{\def\startpage{#1}\def\finishpage{#2}}
\def\published#1{\def\publishdate{#1}}
\def\proposed#1{\def\theproposer{#1}}
\def\seconded#1{\def\theseconders{#1}}
\def\received#1{\def\receiveddate{#1}}

\def\accepted#1{\def\accepteddate{#1}}
\def\asciititle#1{\def\theasciititle{#1}}
\def\covertitle#1{\def\thecovertitle{#1}}

\long\def\asciiabstract#1{\long\def\theasciiabstract{#1}}
\def\asciikeywords#1{\def\theasciikeywords{#1}}

\let\\\par\let\thelognumber\relax
\let\thevolumenumber\relax\let\thepapernumber\relax
\let\thevolumeyear\relax\let\thesamplenumber\relax\let\startpage\relax
\let\finishpage\relax\let\publishdate\relax\let\receiveddate\relax
\let\reviseddate\relax\let\accepteddate\relax\let\theasciititle\relax
\let\thecovertitle\relax\let\theasciiauthors\relax
\let\theasciiabstract\relax\let\theasciikeywords\relax
\let\theasciiemail\relax\let\theshortauthors\relax\let\theshorttitle\relax

\long\def\maketitlep{   

\count0=\startpage

\gt\hfill      
\beginpicture
\setcoordinatesystem units <0.33truein, 0.33truein> point at 2.2 0.9
\setplotsymbol ({$\cal G$})
\plotsymbolspacing=9truept
\circulararc 315 degrees from 0 1 center at 0 0
\setplotsymbol ({$\cal T$})
\circulararc 315 degrees from 1 -1 center at 1 0
\endpicture
\break
{\small\ifx\thesamplenumber\relax 
Volume \else Sample
\fi\thevolumenumber\ (\thevolumeyear)
\startpage--\finishpage\nl
Published: \publishdate}
\vglue 0.5truein plus 0.4fil minus 0.1truein

{\parskip=0pt\leftskip 0pt plus 1fil\def\\{\par\smallskip}{\ifplaintex\large
\else\Large\fi\bf\thetitle}\par\medskip}   

\vglue 0pt plus 0.1fil 

{\parskip=0pt\leftskip 0pt plus 1fil\def\\{\par}{\sc\theauthors}
\par\medskip}

\vglue 0pt plus 0.1fil 

{\small\parskip=0pt\let\newline\\
{\leftskip 0pt plus 1fil\def\\{\par}{\sl\theaddress}\par}
\expandafter\ifx\theemail\relax    
\relax\else\vglue 5pt plus 0.02fil minus 2pt\def\\{\stdspace{\rm 
and}\stdspace} 
\cl{Email:\stdspace\tt\theemail}\fi
\ifx\theurl\relax                  
\relax\else\vglue 5pt plus 0.02fil minus 2pt\def\\{\stdspace{\rm 
and}\stdspace}
\cl{URL:\stdspace\tt\theurl}\fi\par}

\vglue 7pt plus 0.3fil minus 3pt

{\bf Abstract}
\vglue 5pt plus 0.1fil minus 2pt

\theabstract

\vglue 7pt plus 0.3fil minus 3pt

{\bf AMS Classification numbers}\quad Primary:\quad \theprimaryclass

Secondary:\quad \thesecondaryclass

\vglue 5pt plus 0.3fil minus 2pt

{\bf Keywords:}\quad \thekeywords

\vglue 10pt plus 0.5fil minus 5pt

{\small  Proposed: \theproposer\hfill Received: \receiveddate\nl
Seconded: \theseconders\hfill 
\ifx\reviseddate\relax                         
Accepted: \accepteddate                        
\else
Revised: \reviseddate                          
\fi}
\eject
}       

\let\maketitlepage\maketitlep
\let\maketitle\maketitlepage

\font\phead=cmsl9 scaled 950
\font=cmsl9 scaled 1050
\font\pnum=cmbx10 scaled 913
\font\lnum=cmbx10 
\font\pfoot=cmsl9 scaled 950
\font=cmsl9 scaled 1050
\ifplaintex
\headline{\vbox to 0pt{\vskip -4.5mm\line{\small\phead\ifnum
\count0=\startpage ISSN 1364-0380 (on line)
1465-3060 (printed) \hfill {\pnum\folio}\else\ifodd\count0\def\\{ }%
\ifx\theshorttitle\relax\thetitle\else\theshorttitle\fi\hfill{\pnum\folio}
\else\def\\{ and }{\pnum\folio}\hfill\ifx\theshortauthors\relax\theauthors
\else\theshortauthors\fi\fi\fi}\vss}}
\footline{\vbox to 0pt{\vglue 0mm\line{\small\pfoot\ifnum\count0=\startpage
\copyright\ \gtp\hfill\else
\gt, Volume \thevolumenumber\ (\thevolumeyear)\hfill\fi}\vss
}}
\else
\makeatletter
\def\@oddhead{{\small\lhead\ifnum\count0=\startpage ISSN 1364-0380 (on line)
1465-3060 (printed) \hfill {\lnum\number\count0}\else\ifodd\count0
\def\\{ }\ifx\theshorttitle\relax \thetitle \else\theshorttitle\fi\hfill
{\lnum\number\count0}\else\def\\{ and }{\lnum\number\count0}
\hfill\ifx\theshortauthors\relax 
\theauthors\else\theshortauthors\fi\fi\fi}}\def\@evenhead{\@oddhead}
\def\@oddfoot{\small\lfoot\ifnum\count0=\startpage\copyright\ \gtp\hfill\else
\gt, Volume \thevolumenumber\ (\thevolumeyear)\hfill\fi}
\def\@evenfoot{\@oddfoot}
\makeatother
\fi

\newwrite\gtoutfile
\long\gdef\makeheadfile{  
{\def\\{, }\def\s{ }
\immediate\openout\gtoutfile head.xxx
\immediate\write\gtoutfile{Proxy-for: \ifx\theasciiauthors\relax
\theauthors\else\theasciiauthors\fi\s<\ifx\theasciiemail\relax\theemail\else\theasciiemail\fi>}
\immediate\write\gtoutfile{\noexpand\\}
\immediate\write\gtoutfile{Authors: \ifx\theasciiauthors\relax
\theauthors\else\theasciiauthors\fi}
{\def\\{ }\immediate\write\gtoutfile{Title: \ifx\theasciititle\relax
\thetitle\else\theasciititle\fi}}
\immediate\write\gtoutfile{Subj-class: GT or SG or MG etc}
\immediate\write\gtoutfile{MSC-class: \theprimaryclass\ifx\thesecondaryclass\relax\else, \thesecondaryclass\fi}
\immediate\write\gtoutfile{Journal-ref: Geom. Topol. \thevolumenumber
(\thevolumeyear) \startpage-\finishpage}
\immediate\write\gtoutfile{Comments: Published by Geometry and Topology at}
\immediate\write\gtoutfile{\s\s http://www.maths.warwick.ac.uk/gt/GTVol\thevolumenumber/paper\thepapernumber.abs.html}
\immediate\write\gtoutfile{\noexpand\\}
\immediate\write\gtoutfile{}
\ifx\theasciiabstract\relax
\immediate\write\gtoutfile{\theabstract}\else
\immediate\write\gtoutfile{\theasciiabstract}\fi
\immediate\write\gtoutfile{}
\immediate\write\gtoutfile{\noexpand\\}
\immediate\write\gtoutfile{}
\immediate\closeout\gtoutfile}}  

\def\maketitlepage{\maketitlep\makeheadfile}
\let\maketitle\maketitlepage

\lognumber{453}
\received{11 May 2004}
\volumenumber{9}\papernumber{5}\volumeyear{2005}
\pagenumbers{179}{202}   
\published{8 January 2005}
\accepted{27 December 2004}
\proposed{Benson Farb}
\seconded{Jean-Pierre Otal, Walter Neumann}

\usepackage{amssymb, amsmath, graphicx, psfrag}

\psfrag{a}{\small$a$}
\psfrag{b}{\small$b$}
\psfrag{c}{\small$c$}
\psfrag{n=2}{\small$n=2$}
\psfrag{n=1}{\small$n=1$}
\psfrag{n=-1}{\small$n=-1$}

\def\S{Section }

\newtheorem{theorem}{Theorem}[section]

\newtheorem{corollary}[theorem]{Corollary}
\newtheorem{lemma}[theorem]{Lemma}

\theoremstyle{definition}
\newtheorem{definition}[theorem]{Definition}
\newtheorem{example}[theorem]{Example}

\theoremstyle{remark}
\newtheorem{remark}[theorem]{Remark}

\newcommand{\thmref}[1]{Theorem~\ref{#1}}

\newcommand{\secref}[1]{\S\ref{#1}}
\newcommand{\lemref}[1]{Lemma~\ref{#1}}
\newcommand{\corref}[1]{Corollary~\ref{#1}} 
\newcommand{\figref}[1]{Figure~\ref{#1}}
\newcommand{\exref}[1]{Example~\ref{#1}}
\newcommand{\remref}[1]{Remark~\ref{#1}}

\newcommand{\startproof}[1]{\proof[Proof of #1]}
\let\finishproof\endproof

\newcommand{\area}{\mathrm{area}}
\newcommand{\Mod}{\mathrm{Mod}}

\newcommand{\col}{\colon\thinspace}

\newcommand{\bfheading}[1]{\par\medskip{\bf #1}\qua}

\newcommand{\CC}{{\mathcal C}}

\newcommand{\T}{{\mathcal T}}
\newcommand{\TT}{\hat {\mathcal T}}

\newcommand{\F}{\mathcal{PMF}}

\newcommand{\CS}{{\mathcal CS}}

\newcommand{\Teich}{Teichm\"uller }

\newcommand\R{{\mathbb R}}

\newcommand{\ep}{\epsilon}

\newcommand{\nup}{\nu_+}
\newcommand{\nm}{\nu_-}
\newcommand{\mup}{\mu^+}
\newcommand{\mm}{\mu^-}
\newcommand{\npb}{\bar \nu_+}
\newcommand{\nmb}{\bar \nu_-}

\newcommand{\p}{\partial}

\newcommand{\h}{{\bf h}}

\begin{document}

\title{A characterization of short curves\\of a Teichm\"uller geodesic} 
\asciititle{A characterization of short curves of a Teichmueller geodesic} 
\covertitle{A characterization of short curves\\of a 
Teichm\noexpand\"uller geodesic} 
\author{Kasra Rafi}
\address{Department of Mathematics, University of 
Connecticut\\Storrs, CT 06269, USA}
\email{rafi@math.uconn.edu}

\begin{abstract}
We provide a combinatorial condition characterizing curves that are
short along a \Teich geodesic. This condition is closely related to the
condition provided by Minsky for curves in a hyperbolic
3--manifold to be short. We show that short curves in a hyperbolic
manifold homeomorphic to $S \times \R$ are also short in the 
corresponding \Teich geodesic, and we provide examples demonstrating
that the converse is not true.
\end{abstract}

\asciiabstract{%
We provide a combinatorial condition characterizing curves that are
short along a Teichmueller geodesic. This condition is closely related
to the condition provided by Minsky for curves in a hyperbolic
3-manifold to be short. We show that short curves in a hyperbolic
manifold homeomorphic to S x R are also short in the corresponding
Teichmueller geodesic, and we provide examples demonstrating that the
converse is not true.}

\asciikeywords{Teichmueller space, geodesic, short curves, complex of curves, Kleinian group, bounded geometry} 

\maketitlepage

\section{Introduction}  \label{sec:intro}

\noindent
We are interested in studying the behavior of geodesics in \Teich space. 
Let $S$ be a surface of finite type, $\T(S)$ be the \Teich space of $S$ and 
$\TT(S)$ be the Thurston compactification of $\T(S)$. The boundary of 
$\TT(S)$ is the space of projectivized measure foliations, $\F(S)$. A  geodesic 
$g$ in $\T(S)$ can be described by a one-parameter family of quadratic 
differentials $q_t$ on $S$. The horizontal and the vertical foliations of $q_t$ 
define two elements, $\mup$ and $\mm$, in $\F(S)$ which are independent of 
$t$. We give a characterization of short curves of $(g,\mup, \mm)$ based on 
the topological type of $\mup$ and $\mm$ only. Define 
$$l_g(\alpha) = \inf_{\sigma \in g} l_{\sigma}(\alpha)$$
to be the shortest hyperbolic length of $\alpha$ over all hyperbolic metrics 
along the geodesic $g$. For a fixed $\ep > 0$, we say $\alpha$ is {\it short} 
along  $g$ if $l_g(\alpha)\leq \ep$. Also, for a fixed $K > 0$ we say two 
essential arcs in a surface are {\it quasi-parallel} if their intersection number 
is less than $K$. The following is our main theorem:

\begin{theorem} \label{main}
Let $(g,\mup,\mm)$ be a geodesic in $\T(S)$ and $\alpha$ be a simple 
closed curve in $S$. Then $\alpha$ is short along $g$ if and only if, for 
some component $Y$ of $X \setminus \alpha$, the restrictions of  
$\mup$ and $\mm$ to $Y$ have no common quasi-parallels.\footnote{See
\thmref{main-long} and \remref{h=0}.}
\end{theorem}

\begin{corollary}
The set of short curves in $(g,\mup,\mm)$ is independent of the measure 
class given on $\mup$ and $\mm$.
\end{corollary}

This corollary sheds some light on the behavior of geodesics
whose ``endpoints" are not uniquely ergodic. We will
discuss this and other applications of \thmref{main} in
studying the geometry of $\T(S)$ and its stable and unstable 
geodesics in future papers.

The proof of the main theorem requires close attention to the 
geometry of the singular Euclidean metric defined on $S$ by 
a quadratic differential. An essential step is to prove an analogue of 
the collar lemma, which compares the lengths, in quadratic differential 
metric, of intersecting curves, assuming one of them has bounded
hyperbolic length. Let $X$ be a Riemann  surface, $q$ be a quadratic 
differential on $X$ and $\sigma$ be the corresponding hyperbolic metric 
on $X$.

\begin{theorem} \label{q-collar}
For every $L>0$, there exists  $D_L$ such that, if $\alpha$ and $\beta$ 
are two simple closed curves in $X$  intersecting nontrivially with 
$l_{\sigma}(\beta) \leq L$, then
$$  D_L \, l_q(\alpha) \geq l_q(\beta). $$
\end{theorem}
 
\bfheading{Relation with Kleinian groups}
The original motivation for this work was the question of whether the 
universal curve over a \Teich geodesic can be used to model an 
infinite-volume hyperbolic 3--manifold. Here we briefly discuss the 
connection between geodesics in $\T(S)$ and hyperbolic 3--manifolds and
present some difficulties one might encounter in attempting to make such 
a model. In \cite{rees}, Rees recently proposed a solution to these 
difficulties.

Let $M=S \times [0,1]$, and let $\h(M)$ be the space of complete 
hyperbolic structures on the interior of $M$. Elements of $\h(M)$ are 
infinite-volume hyperbolic $3$--manifolds homeomorphic to $S \times \R$. 
The intersection of $N \in \h(M)$ with a neighborhood of a boundary 
component of $M$ is called an {\it end} of $N$. With each end one can 
associate an invariant $\mu$ in an enlargement of $\T(S)$ related to 
Thurston's compactification of the Teichm\"uller space. If an end is 
{\it geometrically finite}, that is, if the conformal structure on 
$S \times \{ t \}$ stabilizes as $t$ approaches infinity, $\mu$ is defined to 
be the limiting conformal structure on $S$, i.e., it is a point in $\T(S)$.
If an end is degenerate, then $\mu$ is an ``unmeasured" foliation on 
$S$.\footnote{An end invariant, in general, is a ``hybrid" of a conformal 
structure and an unmeasured foliation.} The ending lamination theorem 
states that these invariants are sufficient to determine the hyperbolic 
$3$--manifold $N$.

\begin{theorem}[Brock, Canary, Minsky]\label{elc} 
The topological type of $N$ and its end invariants determine 
$N$ up to isometry.
\end{theorem}

The proof of the theorem involves building a model manifold
for the hyperbolic 3--manifold $(N, \mup, \mm) \in \h(M)$ based
on combinatorial information given by $\mup$ and $\mm$ (see \cite{elcI}).
In particular, short curves in $N$ are exactly the boundary components 
of surfaces $Y$ such that projections of $\mup$ and $\mm$ to $Y$ are far 
apart in the complex of curves of $Y$, that is, $d_Y(\mup, \mm)$ is large 
(see \cite{CCII} for definition and discussion). This condition is very closely 
related to the combinatorial condition given in \thmref{main}. Consider the 
geodesic $(g, \mup, \mm)$ in $\T(S)$. (We are abusing the notation 
here, since end invariants of $N$ are unmeasured foliations and
$\mup$ and $\mm$ are in $\F(S)$.) We say $N$ has bounded geometry 
if the injectivity radius of $N$ is positive, and we say $g$ has bounded 
geometry if, for some $\ep > 0$, the injectivity radius of the hyperbolic 
metric at $g(t)$ is larger than $\ep$. Using the technology 
developed in \cite{CCII} and \cite{elcI} we prove:

\begin{theorem} \label{one}
The following are equivalent:
\begin{enumerate}
\item $N$ has bounded geometry,
\item $g$ has bounded geometry.\footnote{This theorem is also a corollary 
of \thmref{elc} and \cite{m-teich}.}
\end{enumerate}
\end{theorem}

\begin{theorem} \label{two}
For every $\ep $ there exists $\ep'$ such that, if the length of
$\alpha$ in $N$ is less than $\ep'$, then
$l_g(\alpha) \leq \ep$.
\end{theorem}

The converse of the last theorem is not true. In
\secref{sec:cc} we will provide a family of counterexamples.

\begin{theorem} \label{counter}
There exist a sequence of manifolds $(N_n,\mup_n,\mm_n)$,
corresponding Teichm\"uller geodesics $(g_n, \mup_n, \mm_n)$ and 
a curve $\gamma$ in $S$ such that $l_{N_n}(\gamma) \geq c>0$
for all $n$, but
$$l_{g_n}(\gamma) \to 0  \qquad \text{as} \qquad  n \to \infty.$$
\end{theorem}

\bfheading{Outline of the paper}
We begin in \secref{sec:curve} by introducing subsurface projections
and reviewing some of their properties. In \secref{sec:quad}
we discuss the geometry of a quadratic differential metric
on a Riemann surface and the behavior of geodesics there. In
\secref{sec:collar} we prove the ``collar lemma" (\thmref{q-collar}).
Theorems establishing a connection between
the combinatorial properties of horizontal and vertical foliations
and the short curves in the quadratic differential metric are
presented in \secref{sec:thm}. The main theorem is proven 
in \secref{sec:proof}, and Theorems \ref{one}, \ref{two}
and \ref{counter} are proven in \secref{sec:cc}.

\bfheading{Notation}
To simplify our presentation we use the notations $O$, $\succ$ and 
$\asymp$, defined as follows: for two functions $f$ and $g$, 
$f \succ g$ means $f \geq c \, g - d$, $f \asymp g$ means 
$c \, f - d \leq g \leq C \, f + D$, and $f = O(g)$ means $ f \leq C g$, 
where $c,d,C$ and $D$ depend on the topology of $S$ only. 

\bfheading{Acknowledgment}
I would like to thank the referee for careful reading and for 
numerous comments and corrections.


\section{Simple closed curves} \label{sec:curve}

\noindent
In this section we study some properties of arcs and curves on
surfaces. Let $S$ be an orientable surface of finite type, excluding the 
sphere and the torus.
By a {\it curve} we mean a non-trivial, non-peripheral, piecewise smooth
simple closed curve in $S$. The free homotopy class of a curve $\alpha$ 
is denoted by $[\alpha]$. By an {\it essential arc} $\omega$ we mean 
a piecewise smooth simple arc, with endpoints on the boundary of $S$, 
that cannot be pushed to the boundary of $S$. In case $S$ is not an annulus,
$[\omega]$ represents the homotopy class of $\omega$ relative to the boundary
of $S$. When $S$ is an annulus, $[\omega]$ is defined to be  the homotopy 
class of $\omega$ relative to the endpoints of $\omega$.

Define $\CC(S)$ to be the set of all homotopy classes of
curves and essential arcs on the surface $S$. A curve system  
$\Gamma = \{ \gamma_1, \ldots, \gamma_n \}$
is a non-empty set of curves and essential arcs in $S$ that are pairwise
disjoint from each other. The homotopy class of $\Gamma$ is the set of
homotopy classes of elements of $\Gamma$,
$[ \Gamma ] =  \{ [\gamma_1], \ldots, [\gamma_n] \}$. 
Let $\CS(S)$ be the space of all classes of curve systems on $S$. 
For two curve systems $\Gamma_1$ and $\Gamma_2$ we define
the intersection number between $\Gamma_1$ and $\Gamma_2$, $i_S(\Gamma_1, \Gamma_2)$, to be the minimum geometric 
intersection number between an element of $\Gamma_1$ and 
an element of $\Gamma_2$. This depends only on the homotopy classes 
of $\Gamma_1$ and $\Gamma_2$. If $S$ is an annulus, $\Gamma_1$
includes the core curve of $S$ and $\Gamma_2$ includes any
transverse curve, then we define their intersection number to be
infinity.

A subsurface $Y$ of $S$ is a connected open subset of $S$ 
with piecewise smooth boundary. Let $[Y]$ denote the homotopy class 
of $Y$, that is, the set of all connected open subsets of $S$ with piecewise
smooth boundary that are homotopically equivalent to $Y$. A choice of
base point in $Y$ identifies $\pi_1(Y)$ with a subgroup of $\pi_1(S)$.
In fact, this subgroup is well defined up to conjugation for $[Y]$.

Let $\mu$ be a curve system, a lamination or a singular-foliation on $S$.
Here we define the {\it projection} of $\mu$ to the subsurface $Y$. Let
$$f \col \bar S \to S$$
be a covering of $S$ such that $f_*(\pi_1(\bar S))$ is conjugate to
$\pi_1(Y)$. We call this a $Y$--cover of $S$. Since $S$ admits a hyperbolic 
metric, $\bar S$ has a well defined boundary at infinity. Let $\bar \mu$ be the 
lift of $\mu$ to $\bar S$. Components of $\bar \mu$ that are curves or 
essential arcs on $\bar S$, if any, form a curve system in $\bar S$. This set, 
if not empty, defines an element in $\CS(\bar S)$. Since $\bar S$ is 
homeomorphic to $Y$, $\CS(\bar S)$ is identified with $\CS(Y)$ through 
this homeomorphism.  We call the corresponding element of $\CS(Y)$ the 
projection of $\mu$ to $Y$ and will denote it by $\mu_Y$. If there are 
no essential arcs or curves in $\bar \mu$, $\mu_Y$ is the empty set; 
otherwise we say that $\mu$ intersects $Y$ non-trivially. Note that this 
projection depends on the homotopy class of $\mu $ only.
 
\bfheading{Subsurface intersection}
Let $\mu$ and  $\nu$ be two curve systems, laminations or singular foliations 
on $S$ that intersect a subsurface $Y$ non-trivially. We define the
{\it $Y$--intersection} of $\mu$ and $\nu$ to be the intersection number in 
$Y$ between the projections $\mu_Y$ and $\nu_Y$ and denote it by
$$i_Y(\mu,\nu) = i_Y(\mu_Y, \nu_Y).$$


\section{Quadratic differentials}
\label{sec:quad}

\noindent
In this chapter we review the geometry of quadratic differentials. 
For more detailed discussion, see \cite{m-thesis} and \cite{strebel}.

Let $\hat X$ be a compact Riemann surface  and let $q \not \equiv 0$ be a
meromorphic quadratic differential on $\hat X$. We assume that the set of 
critical points of $q$ is discrete and that $q$ has finite critical points 
only (i.e., critical points of $q$ are either zeros or poles of order 1). 
The surface $\hat X$ punctured at poles $P_i$ of $q$ is denoted by $X$.
For any such $q$, $|q|$ defines a singular Euclidean 
metric in $X$. We denote the length of a piecewise smooth curve by $l_q(\gamma)$ and the area of a subset $Y$ of $X$ by $\area_q(Y)$. 
We also assume that $\area_q(X)=1$. The metric defined on $X$
is not complete since the finite poles are finite distance from interior 
points of $X$. Let $\tilde X$ be the completion of the universal cover of 
$X$. In a neighborhood of every regular point $P$ of $q$ 
we can introduce a local parameter $w$, in terms of which the representation 
of $q$ is identically equal to one. This parameter, given by the integral
$$w=Q(z)= \int \sqrt{q(z)} \, dz,$$
is uniquely determined up to a transformation $w \to \pm w + \text{const}$,
and it will be called the natural parameter near $P$.

\bfheading{Horizontal and vertical foliations}
A {\it straight arc} is a smooth curve $\gamma$ such that, for every 
point $P$ in $\gamma$, the image of $\gamma$ under the natural
parameter near $P$ is a straight line. Since $w$ is ``unique", the angle of 
$\gamma$ is a well-defined number $0\leq \theta < \pi$.
A maximal such straight arc is called a {\it $\theta$--trajectory}. 
A straight arc connecting two critical points in $(X,q)$ is called a
{\it straight segment}. The $\theta$--trajectory passing through a regular 
point of $q$ is unique.  A trajectory is called a critical trajectory if at least 
one of its ends tends to a critical point. There are only finitely many critical 
$\theta$--trajectories.

\begin{figure}[ht!]
\begin{center}
\includegraphics[width=3.5in]{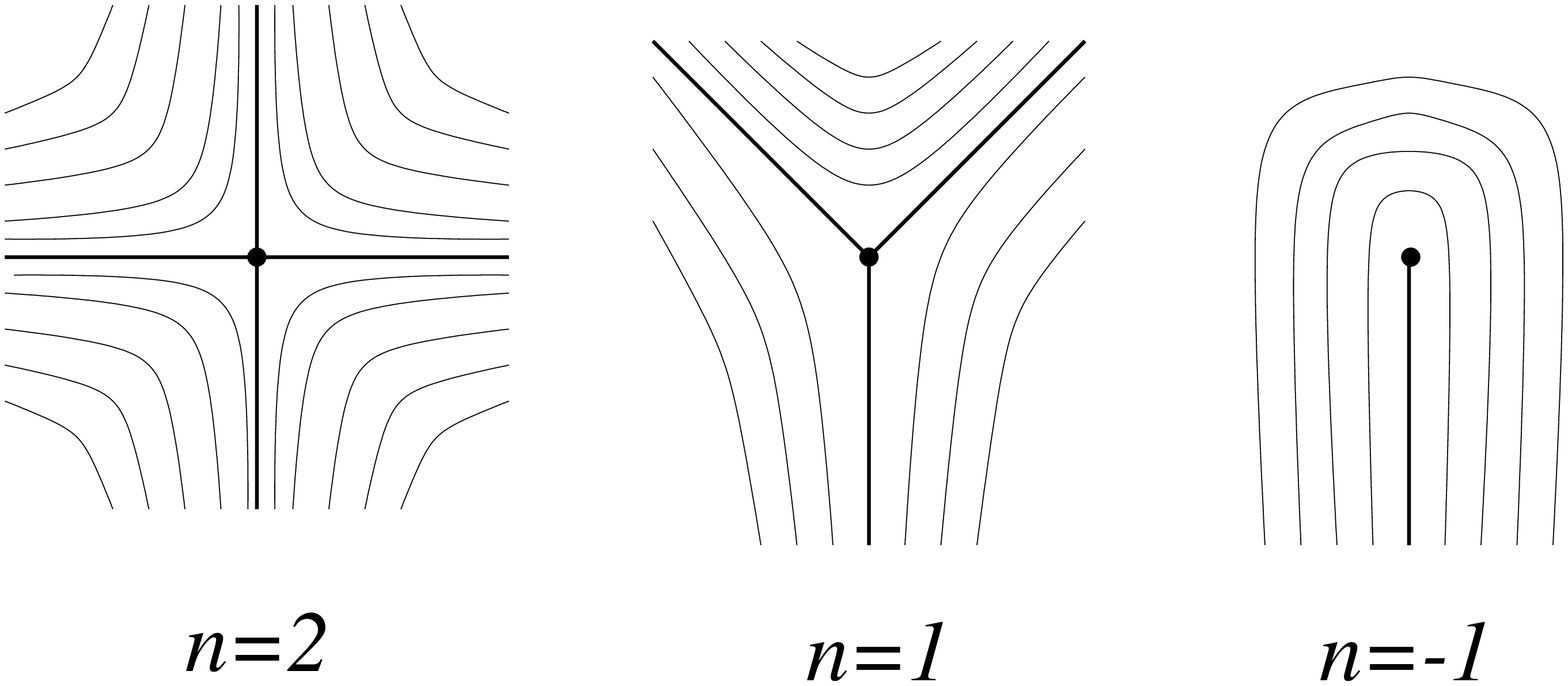}
\caption{Local pictures of a foliation by trajectories} 
\label{fig:sign}
\end{center}
\end{figure}

For each $\theta$, the $\theta$--trajectories foliate the set of
regular points in $X$. For $n=-1,1,2$, a local neighborhood
of a finite critical point of degree $n$ is shown in \figref{fig:sign}.

The foliation corresponding to $\theta = 0$ (respectively, $\theta = \pi /2$), 
is called the {\it horizontal foliation} (the {\it vertical foliation}) and is denoted 
by $\nup$ ($\nm$).  Let $\gamma$ be any closed curve or arc. 
Chop $\gamma$ to coordinate patch sized segments $\gamma_i$.  
For each such segment, we define the horizontal length to be
$$h_q(\gamma_i) =\int_0^1 
| \sqrt{q(\gamma_i(t))} | \, \cdot \, | 
\text{Re} \{ \gamma_i'(t) \} | \, dt.$$ 
Then, the horizontal length of $\gamma$, $h_q(\gamma)$, is defined as 
the sum of horizontal lengths of the $\gamma_i$.
The vertical length of $\gamma$ is defined similarly and is denoted by 
$v_q(\gamma)$. Let $\gamma$ be an arc that is transverse to the horizontal 
foliation. The vertical length of $\gamma$ does not 
change after a homotopic change through arcs that are transverse to the 
horizontal foliation. Therefore the horizontal length defines a transverse measure 
on $\nm$; similarly, the vertical length defines a transverse measure on $\nup$, 
making $\nup$ and $\nm$ into measured foliations.

\bfheading{Geodesic representatives}
We would like to define the {\it $q$--geodesic representative} of a curve 
$\gamma$ to be the curve in the free homotopy class of $\gamma$ with 
the shortest $q$--length. If $X$ has no punctures, the $q$--metric is complete 
and the geodesic representative always exists.  In general, the $q$--metric
on $X$ is not complete; however, for any lift $\tilde \gamma$ of $\gamma$,
there exists a unique geodesic $\tilde \gamma_q$ in $\tilde X$ 
that has the same endpoints at infinity as $\tilde \gamma$.
The projection of $\tilde \gamma_q$ to $\hat X$ is denoted by $\gamma_q$ 
and is called the $q$--geodesic representative of $\gamma$. 
The curve $\gamma_q$ is composed of straight segments which meet 
at critical points of $X$, making an angle of at least $\pi$ on either side.
If $\gamma_q$ ``passes through" a pole of $\hat X$, it makes at least a full 
turn around the pole.

Unfortunately, one cannot recover the homotopy class of $\gamma$
from the set $\gamma_q$. Therefore, we should think of $\gamma_q$
as a parametrized curve. Identify $\tilde \gamma_q$ with 
the real line as a metric space. The deck transformations fixing 
$\tilde \gamma_q$ define a $\mathbb{Z}$ action on it. The quotient is a circle 
of length $l_q(\gamma)$. The projection from $\tilde X$ to $\hat X$
defines a map 
$$f_{\gamma} \col S^1 \to \gamma_q,$$
which we call the natural parametrization of $\gamma_q$.

\begin{theorem}[see \cite{m-thesis}, \cite{strebel}]\label{q-geo}
For $(X,q)$ as above,
\begin{enumerate}
\item Let $\gamma$ be an arc joining two given points $x$ and $y$ of $X$;
then the $q$--geodesic representative of $\gamma$ exists and is unique.
\item Let $\gamma$ be an essential loop on $X$; then the 
$q$--geodesic representative of $\gamma$ exists and is unique, 
except for the case where it is one of the continuous family of closed 
Euclidean geodesics in a flat annulus (see page~\pageref{flat} for definition
of a flat annulus). 
\end{enumerate}
\end{theorem}

A geodesic curve or arc $\gamma$ is called {\it balanced} if 
$v_q(\gamma)=h_q(\gamma)$. If $v_q(\gamma) \geq h_q(\gamma)$
we say $\gamma$ is {\it mostly vertical} and if
$v_q(\gamma) \leq h_q(\gamma)$ we say $\gamma$ is mostly
horizontal. Consider $\gamma$ as a union of straight segments 
$\gamma_i$. Let $l_i$, $h_i$ and $v_i$ be the $q$--length, the 
horizontal length and the vertical length of $\gamma_i$ respectively. 
We have $l_i^2=h_i^2 + v_i^2$. The Cauchy--Schwarz inequality implies
\begin{eqnarray*}
l_q(\gamma)^2 & =     & \left( \sum l_i \right)^2  \leq 2 \sum l_i^2 
                                      = 2 \sum h_i^2 + 2 \sum v_i^2 \\
                         & \leq & 2 \left( \sum h_i \right)^2 + 2\left( \sum v_i \right)^2
                                      = 2 \, (h_q(\gamma)^2 + v_q(\gamma)^2).
\end{eqnarray*}
Therefore, if $\gamma$ is mostly vertical (respectively, mostly horizontal),
then
\begin{equation} \label{mostly}
2 \, v_q(\gamma) \geq l_q(\gamma) \quad 
\big( {\rm respectively,} \quad 2 \, h_q(\gamma) \geq l_q(\gamma) \big).
\end{equation}

\bfheading{Curvature of a curve in $X$}
Let $\gamma$ be a boundary component of $Y$.
The curvature of $\gamma$ with respect to $Y$, $\kappa_Y(\gamma)$, is well 
defined as a measure with atoms at the corners (that is, $(\pi - \theta)$, 
where $\theta$ is the interior angle at the corner). We choose the sign to 
be positive when the acceleration vector points into $Y$. If $\gamma$ is 
curved positively (or negatively) with respect to $Y$ at every point, we say 
it is {\it monotonically curved} with respect to $Y$.

\begin{remark}
The curvature of any closed curve is a multiple of $\pi$(see \cite{m-thesis}). 
\end{remark}

\begin{theorem}[Gauss--Bonnet]\label{GB}
Let $Y$ be a subsurface of $\hat X$ with boundary components 
$\gamma_1, \ldots , \gamma_n$; also let $P_1, \ldots, P_k$
be the critical points in $Y$.  Then we have
$$  \sum_{i=1}^n \kappa_Y(\gamma_i) =  
                 2\pi \chi(Y)+ \pi  \left( \sum_{i=1}^k \deg(P_i) \right),$$
where $\deg(P_i)$ is the order of the critical point $P_i$.
\end{theorem}

\bfheading{Regular and primitive annuli in $X$}
In general, it is difficult to estimate the modulus of an annulus
in $X$. Here we discuss estimates for moduli of certain kinds of
annuli, which will be useful in estimating the hyperbolic
length of a closed curve (see \cite[Section 4]{m-thesis}).

\begin{definition}
Let $A$ be an open annulus in $X$ with boundaries $\gamma_0$ and
$\gamma_1$. Suppose both boundaries are monotonically curved
with respect to $A$. Further, suppose that the boundaries are
equidistant from each other, and $\kappa_A(\gamma_0) \leq 0$.
We call $A$ a {\it regular} annulus. If $\kappa_A(\gamma_0) < 0$,
we call $A$ {\it expanding}  and say that $\gamma_0$ is the inner 
boundary and $\gamma_1$ is the outer boundary. If the interior of
$A$ contains no zeroes, we say $A$ is a {\it primitive} annulus, and we write
$\kappa(A)=-\kappa_A(\gamma_0)$. When $\kappa(A)=0$, $A$ is a {\it flat}
\label{flat} annulus and is foliated by closed Euclidean geodesics 
homotopic to the boundaries.
\end{definition}

\begin{remark}
Later in this section we will discuss certain covers of $X$ which are
not compact. There, a regular annulus can be infinite, having only one 
boundary component, the inner boundary.
\end{remark}

\begin{lemma} \label{p-dist}
Let $A$ be a primitive annulus, $\gamma_0$ and $\gamma_1$ be the inner 
and the outer boundary of $A$ respectively and $d$ be the distance between 
the boundaries of $A$. Then
\begin{equation*}
\left\{ \begin{array}{ll} \displaystyle
\kappa \, \Mod(A) \asymp \log \left( \frac{d}{ l_q(\gamma_0)} \right) &  
                                        \text{if} \: \: \kappa(A) > 0, \\ & \\
\Mod(A) \, l_q(\gamma_0) = d(\gamma_0, \gamma_1)  & \text{if} \: \: \kappa(A) = 0.
\end{array} \right. 
\end{equation*}
\end{lemma}

\begin{proof}
The statement is clear for a flat annulus. In case $\kappa(A) > 0$, Minsky
has shown (\cite[Theorem 4.5]{m-thesis}) that 
$$\Mod(A) \asymp \frac{1}{\kappa(A)} 
\log  \left( \frac{l_q(\gamma_1)}{l_q(\gamma_0)} \right).$$
We have to show that $l_q(\gamma_1)/l_q(\gamma_0) \asymp d/ l_q(\gamma_0)$.
For $0 \leq t \leq d$, define  $\alpha_t$ to be the equidistant curve from 
$\alpha_0$ with $d(\alpha_0, \alpha_t) = t$. Let $A_t$ be the annulus 
bounded by $\alpha_0$ and $\alpha_t$. The Gauss--Bonnet theorem implies 
that $\kappa_{A_t}(\alpha_t)+\kappa_{A_t}(\alpha_0) = 0$. Therefore,
$\kappa_{A_t}(\alpha_t) = \kappa$. The $q$--length of $\alpha_t$
is a differentiable function of $t$ almost everywhere, and
$$\frac{d}{dt}l_q(\alpha_t) = \kappa_{A_t}(\alpha_t),$$
that is, $l_q(\alpha_t) = l_q(\alpha_0) + \kappa \, t$. 
Therefore, $l_q(\gamma_1)= l_q(\gamma_0) + \kappa \, d$.
This finishes the proof.
\end{proof}

We also recall the following theorem of Minsky. This is
a consequence of Theorem 4.6 and Theorem 4.5 in \cite{m-thesis}.
(The inequality $\Mod(A) \leq m_0$ in  \cite[Theorem 4.6]{m-thesis}
is a typo and should read $\Mod(A) \geq m_0$.)

\begin{theorem}[Minsky \cite{m-thesis}]\label{prim-ann}
If $A \subset X$ is any homotopically nontrivial annulus, 
then $A$ contains a primitive annulus $B$ such that 
$$ \Mod(B) \asymp \Mod(A). $$
\end{theorem}

\bfheading{Subsurfaces of $X$}
Given a subsurface $Y$,  it is desirable to consider a representative of
homotopy class $Y$ with $q$--geodesic boundary. The na\"ive approach 
(which would work for the hyperbolic metric) would be to pick a geodesic
representative for each boundary component of $Y$ and take the 
complementary component of these curves in $X$ that is in the homotopy 
class of $Y$. However, the geodesic representatives 
of boundary components of $Y$ might be ``tangent" to each other, that is,
they might share a common geodesic segment. Also, a geodesic representative 
of a simple closed curve could have self-tangencies. As a result, it is possible 
that none of the complementary components are in the homotopy class of $Y$.

Here we describe a representative of $Y$ in the $Y$--cover of
$X$ which, as we shall see, will be more convenient to work with.
Let $\bar X$ be the $Y$--cover of $X$; for each boundary
component $\alpha$ of $Y$,  let $\bar X_\alpha$ be the annular cover
of $X$ corresponding to $\alpha$; and let 
$$p\col \bar X \to X, \qquad \text{and} \qquad 
    p_\alpha \col  \bar X_\alpha \to \bar X $$ 
be the corresponding covering maps. Consider the lift $\hat Y$ of $Y$ to 
$\bar X$ such that $p\co \hat Y \to Y$ is a homeomorphism. The complement 
of $\hat Y$ in $\bar X$ is a union of annuli. We denote the annulus corresponding
to a boundary component $\alpha$ by $B_\alpha$.

The geodesic representative of the core curve of $\bar X_\alpha$
is an embedded two-sided curve dividing $\bar X_\alpha$ into two annuli.
The lift of $B_\alpha$ to $\bar X_\alpha$ shares an end with one of
these annuli, which we denote by $\bar A_\alpha$. The interior of $\bar A_\alpha$ 
is mapped homeomorphically into $\bar X$ by $p_\alpha$.
Otherwise, there would be a disk in $\bar X$ whose boundary
consists of two geodesic arcs. This is  a contradiction, because
the geodesic representation of an arc is unique.  Using the same argument
we also have that the annuli corresponding to different boundary components 
of $\hat Y$ are disjoint.  We denote $p_\alpha(\bar A_\alpha)$ by $A_\alpha$.

If $Y$ is not an annulus 
we define $\bar Y$\label{Ybar} to be the complement in $\bar X$ of all such
annuli $A_\alpha$. Otherwise let $\bar Y$ be the flat annulus containing 
all geodesic representatives of $\alpha$. The surface $\bar Y$ is 
closed and connected, and the restriction of $p$ to the interior of  $\bar Y$ 
is one-to-one. Curves intersecting $Y$ essentially lift to arcs that have a 
compact restriction to $\bar Y$ and wander off to infinity from both 
ends into expanding annuli without ever coming back. That allows us
to talk about the restriction of a curve to $\bar Y$ while keeping track of
the homotopy type of the intersecting curve in the degenerate cases,
where boundaries of $\bar Y$ are tangent to each other.
As a consequence of the Gauss--Bonnet theorem in this setting we have:

\begin{lemma} \label{essential}
For two infinite geodesic arcs $\mu$ and $\nu$ in $\bar X$, the geometric intersection number between $\mu$ and $\nu$ is equal to the intersection
number of the restrictions of $\mu$ and $\nu$ to $\bar Y$ up to an 
additive constant (in fact, the difference is at most 2).
\end{lemma}


\section{Collar lemma for a quadratic differential metric on a surface}
\label{sec:collar}
  
\noindent
Let $X$ and $q$ be as before and $\sigma$ be the corresponding 
hyperbolic metric on $X$. In this section we prove \thmref{q-collar}.

\medskip{\bf \thmref{q-collar}}\qua{\sl
For every $L>0$, there exists  $D_L$ such that if $\alpha$ and $\beta$ 
are two simple closed curves in $X$  intersecting non-trivially with 
$l_{\sigma}(\beta) \leq L$, then
$$  D_L \, l_q(\alpha) \geq l_q(\beta). $$ 
}

Let $B_{\beta}$ be the $\delta$--neighborhood of the geodesic 
representative of $\beta$ in the hyperbolic metric of $X$, where
$\delta$ is a largest number such that this neighborhood is
an annulus, and $m$ be the modulus of $B_\beta$. We have
(see \cite{maskit})
\begin{equation}
\log (m) \succ - L.
\end{equation}
Also let $\ep$ be such that $B_{\beta}$ is disjoint from the $\ep$--thin 
part of the hyperbolic metric on $X$. We can choose $\ep$ such that
$\ep<L$ and $e^L \asymp \frac 1 \ep$.  (The reader should keep in mind
that the theorem is of interest for $L$ large.  The condition $\ep<L$ is
included here simply to ensure that the proof works for small values of $L$ 
as well.) Now pick $M \asymp \frac 1 \ep$ such that, 
if $A$ is an annulus with modulus greater than or equal to $M$, then the 
length of a core curve in $A$ in the complete hyperbolic metric on 
$A$ is less than $\ep$. The Schwarz lemma implies that the $\sigma$--length 
of this curve is also less than $\ep$ and therefore $B_{\beta}$ does not 
intersect the core of $A$. We need the following technical lemma.

\begin{lemma}
Let $\beta$ and $\alpha$ be as above. Then there exists a subsurface $Z$
containing $B_{\beta}$ such that
$$\area_q(Z) \leq R_L \, l_q(\alpha_q)^2$$
where $R_L$ is a constant depending on $L$ and the topology of $X$.
\end{lemma}

\begin{proof}
The strategy is to take $Z$ to be the smallest neighborhood of $\alpha_q$ 
that contains a primitive annulus of modulus $M$ in a neighborhood of
each of its boundary components. Define $N_r$ to be the open 
$r$--neighborhood of $\alpha_q$ in $X$, and let $Z_r$ be the union of $N_r$ 
with all the components of $X-N_r$ that are disks or punctured disks. 
Let $K= - 2 \pi \chi(X) $ and $\kappa_r$ be the sum of the curvatures
of the boundary components of $Z_r$ with respect to $Z_r$. The 
Gauss--Bonnet theorem implies\footnote{Use $\chi(Z_r) \leq 0$ and 
$\sum_{i=1}^k \deg(P_i) \leq -2 \chi(X)$.} 
$\kappa_r \leq  K $.

Observe that $l_q(\p Z_r)$ and $\area_q(Z_r)$ are differentiable
functions of $r$ almost everywhere and are continuous 
except at finitely many points, where we add a disk or a 
punctured disk. For differentiable points we have  
$\frac{d}{dr} l_q(\p Z_r) = \kappa_r$.
At any $r_u$ where we add a disk or a punctured disk $D_u$, the value 
of $l_q(\p Z_r)$ decreases by the $q$--length of the boundary of $D_u$,
which we call $c_u$. Let $I_r$ be the set of all indexes $u$ such that 
$r_u \leq r$. We have the following equation:
\begin{equation}
l_q(\p Z_r) -l_q(\p Z_0) = \int_0^r \kappa_{\rho} \, d \rho - \sum_{u \in I_r} c_u.
\end{equation}
But $l_q(\p Z_0)= 2 l_q(\alpha_q)$ and $\kappa_{\rho} \leq K$. Therefore,
\begin{equation}  \label{lqr}
l_q(\p Z_r) \leq K \, r + 2 l_q(\alpha_q)
\end{equation}
and
\begin{equation}  \label{cu}
\sum_{ u \in I_r } c_u \leq K \, r + 2l_q(\alpha_q).
\end{equation}
To find an upper bound for $\area_q(Z_r)$, we observe that
$$\frac{d}{dr} \area_q(Z_r) = l_q(\p Z_r).$$
At any $r_u$ where we add a disk or a punctured disk,  
$\area_q(Z_r)$ increases by $O(c_u^2)$.\footnote{This is the 
isoperimetric inequality. For the punctured disk, where 
there is a point of positive curvature, consider the double-cover and add 
the missing point.} Therefore,
\begin{equation*}
\area_q(Z_r) -\area_q(Z_0) = \int_0^r l_q(\p Z_{\rho}) \, d\rho + \sum_{u \in I_r} O(c_u^2). 
\end{equation*}
But $\area_q(Z_0)=0$. Equation~\eqref{lqr} and the Cauchy--Schwarz 
inequality imply
$$\area_q(Z_r) \leq \int_0^r \big( K  \rho  \,+ 2l_q(\alpha_q) \big) 
\, d\rho + 
\big[ \sum_{u \in I_{r}} O(c_u) \big]^2 . $$
Now use \eqref{cu} and take the integral to obtain 
\begin{equation}
\label{area}
\area_q(Z_r) \leq \frac{Kr^2}{2} + 2 l_q(\alpha_q) r + 
O\left( \big( K\, r+2l_q(\alpha_q) \big)^2 \right).
\end{equation}

Let $r_0=0$ and $0 < r_1 <  \ldots < r_s$ be points where the topology
of $Z_r$ changes or a critical point of $X$ is added to $Z_r$.
The set $Z_{r_{t+1}} \setminus Z_{r_t}$ is a union of primitive annuli.
Let $t$ be the smallest index such that the moduli of all these primitive 
annuli are greater than $M$. If this never happens, let $t=s$. 

For $j < t$ there exists a primitive annulus $B$ with modulus
less than $M$ connecting a boundary component of $Z_{r_{j+1}}$ to a 
boundary component of $Z_{r_j}$. The distance between the
boundaries of $B$ is $r_{j+1} - r_j$ and the length of the shorter boundary component of $B$ is at most $K r_j + 2 l_q(\alpha_q)$ 
(Equation~\eqref{lqr}). Using \lemref{p-dist}, there exists a $c > 0$ so that
$$ r_{j+1} - r_j \leq c^{K M} \big( K r_j +2l_q(\alpha_q) \big). $$
Let $P=c^{K M} K + 1$ and $Q= 2 \, c^{K M} $. We have
\begin{equation} \label{diff1}
r_{j+1} \leq P \, r_j + Q \, l_q(\alpha_q).
\end{equation}
Now pick $\hat r > r_t$ to be the smallest number such that modulus
of every annulus in $Z_{\hat r} \setminus Z_{r_t}$ is larger than or equal 
to $M$. This is possible by the choice of $t$. If $Z_{r_t}$ equals $X$,
then $\hat r = r_t$, otherwise, by choosing $\hat r$ to be the smallest 
such number, we guarantee the existence of a primitive annulus
of modulus exactly $M$. Then the above argument shows
\begin{equation} \label{diff2}
\hat r \leq P \, r_t + Q \, l_q(\alpha_q) .
\end{equation}
Combining \eqref{diff1} for $0 \leq j < t$ and \eqref{diff2} we get
\begin{equation}
\hat r \leq P^{t +1} r_0 + \sum_{i=0}^t P^i \, Q \, l_q(\alpha_q) . 
\end{equation}
But $r_0=0$. Therefore
$$ \hat r \leq \frac{(P^{t+1} -1)Q}{P-1} \, l_q(\alpha_q).$$
Let $Z=Z_{\hat r}$. Combining this inequality with \eqref{area}
we have, for an appropriate $R_L$, 
$$\area_q(Z) \leq R_L \, l_q(\alpha_q)^2.$$
Note that $R_L$ depends on $L$ or the topology of $X$ only. 
We have to show that $Z$ contains $B_{\beta}$. Let 
$A_1, \ldots , A_k $ be the annuli in $Z_{\hat r} \setminus Z_{r_t}$.
Since $\beta$ intersects $\alpha$
essentially, $B_{\beta}$ has to intersect $Z$ essentially. But the modulus 
of any $A_i$ is larger than $M$; therefore $B_{\beta}$ is disjoint 
from a core curve in each $A_i$. That is, $B_{\beta}$ is contained in $Z$.
\end{proof}

\startproof{\thmref{q-collar}}
Let $\lambda$ be the length of a shortest representative of curve $[\beta]$ 
in $B_{\beta}$ equipped with the $q$--metric inherited from $X$. By definition 
of modulus we have:
\begin{equation}
\frac{\lambda^2}{\area_q(B_{\beta})} \leq \frac1m .
\end{equation}
But $l_q(\beta_q) \leq \lambda$ and $\area_q(B_{\beta}) \leq \area_q(Z)$.
Therefore,
$$\frac{l_q^2(\beta_q)}{\area_q(Z)} \leq \frac{\lambda^2}{\area_q(B_{\beta})}. $$
That is,
$$ l_q^2(\beta_q) \leq \frac{\area_q(Z)}{m}
\leq \frac{R_L l^2_q(\alpha_q)}{m}. $$
Setting $D_L= \sqrt{ \frac{R_L}m }$ finishes the proof of the theorem.
\finishproof

\begin{remark} \label{DL-estimate}
Note that $\log (R_L) \asymp M \asymp e^L$ and 
$\log(m) \asymp -L$. Therefore 
$$\log(D_L) \asymp e^L + L \asymp e^L.$$
\end{remark}
\begin{remark}
The above estimate for $D_L$ is sharp. It is easy to produce examples 
where 
$$\frac{l_q(\beta)}{l_q(\alpha)} \geq e^{e^L}.$$
\end{remark}


\section{Short curves of a quadratic differential metric} \label{sec:thm}

\noindent
Let $X$, $q$, $\nup$ and $\nm$ be as before, $Y$ be a subsurface of $X$ 
and $\bar X$, and $\bar Y$ be as in \secref{sec:quad}. Let  $\bar q$, 
$\npb$ and $\nmb$ be lifts of $q$, $\nup$ and $\nm$ to $\bar X$,
respectively. Fix a positive constant $\delta$.

\begin{definition} 
Let $f \col [0,1] \times (0,1) \to \bar X$ be an embedding with the following 
properties:
\begin{enumerate}
\item For $t \in (0,1)$, $f$ maps $[0,1] \times \{ t \}$ 
into a horizontal leaf of $q$.
\item The intersection of image of $(0,1) \times (0,1)$ under $f$ with $\bar Y$
is non-empty.
\item For $s= 0, 1$, $f$ maps $\{ s \} \times (0,1)$ into 
an arc that is $\delta$--equidistance from $\p \bar Y$ and that is inside one 
of the expanding annuli in $\bar X \setminus \bar Y$.
\end{enumerate}
We call $R=f([0,1] \times (0,1))$ a {\it horizontal strip.} We say a 
horizontal strip is {\it maximal} if it is not a proper subset of 
any other horizontal strips. The width $w(R)$ of a horizontal strip $R$ 
is defined to be the vertical length of a transverse arc to the strip, that is,
$$w(R) = v_q \big( f(\{ t \} \times (0,1)) \big).$$
{\it Vertical strips}, {\it maximal vertical strips} and the width of 
a vertical strip are defined similarly. 
\end{definition}

\begin{lemma} \label{max-strip}
Any two maximal horizontal (respectively, vertical) strips in $\bar Y$
either have disjoint interiors or are identical. Furthermore, 
there are only finitely many distinct maximal strips, and the union of 
all maximal strips covers the boundary of $\bar Y$ except for 
finitely many points. \end{lemma}

\begin{proof}
If the interiors of two horizontal strips intersect, 
their union is also a horizontal strip. But if they are both maximal, 
the union cannot be larger than either of them; therefore they are identical. 
For a point $P$ in the boundary of $\bar Y$, consider the horizontal leaf 
passing through $P$ restricted to a $\delta$--neighborhood of $\bar Y$. 
If this horizontal arc contains no critical points of $q$, then a neighborhood 
of this leaf is also free of critical points, which implies that $P$ is an interior 
point of some horizontal strip. This proves the last assertion.
To see that there are only finitely many strips, note that every maximal strip 
has a critical point on its boundary (otherwise one could extend 
it to a larger strip), and each critical point of degree $n$ can appear on the boundary of at most $n+2$ maximal strips. But the sum of the degrees of 
the critical points in $X$ (and therefore in $\bar Y$) is finite, and so is the 
number of maximal strips. The proof for vertical strips is similar. 
\end{proof}

 \begin{corollary} \label{strip}
Let $R_1, \ldots, R_p$ be the set of maximal horizontal  strips. 
Then the total vertical  length of the boundary of $\bar Y$ 
equals twice the sum of the widths of  $R_i$:
\begin{equation}
v_{\bar q}(\p \bar Y) = 2 \sum_{i=1}^p w(R_i). 
\end{equation}
A similar statement holds for the set of maximal vertical strips. 
\end{corollary}

\begin{example} \label{example}
Let $Y$ be a twice-punctured torus and and $\bar X$ be the following
metric space. Consider a square torus and two copies of the Euclidean 
plane. Cut a vertical slit of the same size (say, of length $l$) at the middle 
of each of them and glue the sides of the slits as shown in \figref{fig:X}. 

\begin{figure}[ht!]
\begin{center}
\includegraphics[width=3.5in]{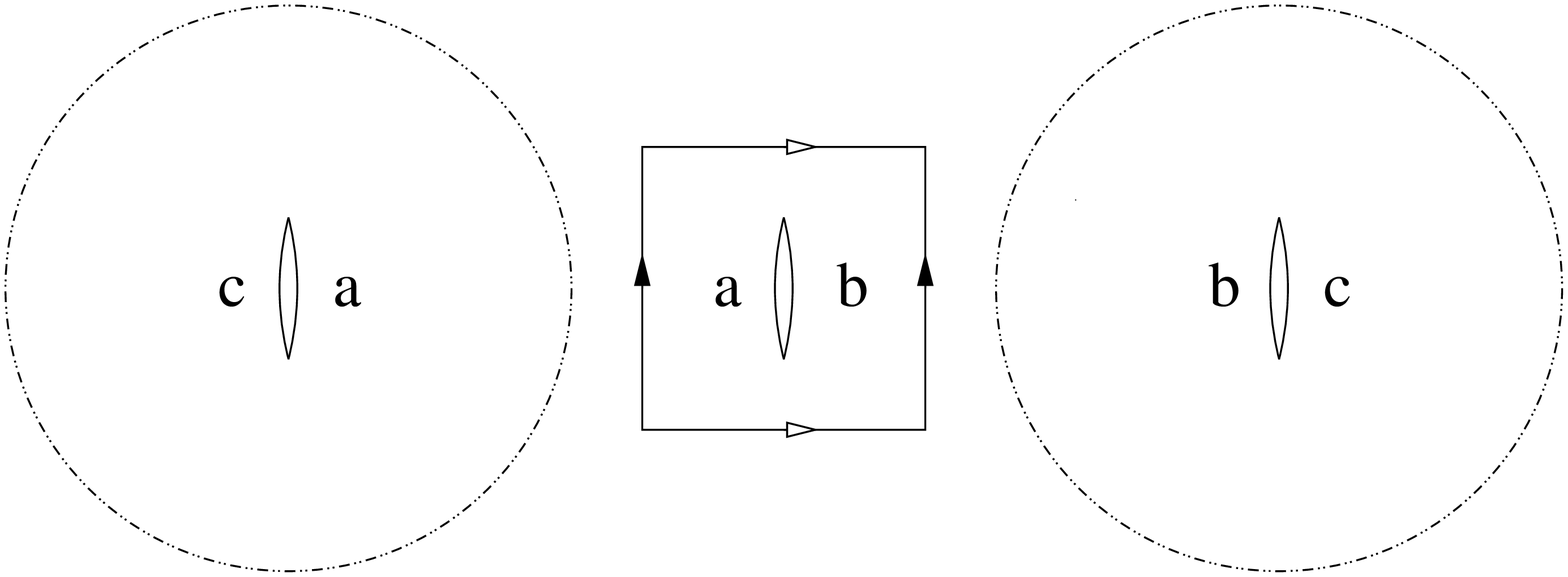}
\caption{Space $\bar X$} 
\label{fig:X}
\end{center}
\end{figure}
This defines a singular Euclidean metric on $\bar X$ with
two singular points of angle $6 \pi$ each. In this case $\bar Y$
is the union of the square torus with the arc $c$. Let the horizontal 
foliation be the foliation by horizontal lines. In \figref{fig:strips} both maximal 
horizontal strips are shown, each  with different shading. 
\begin{figure}[ht!]
\begin{center}
\includegraphics[width=3.5in]{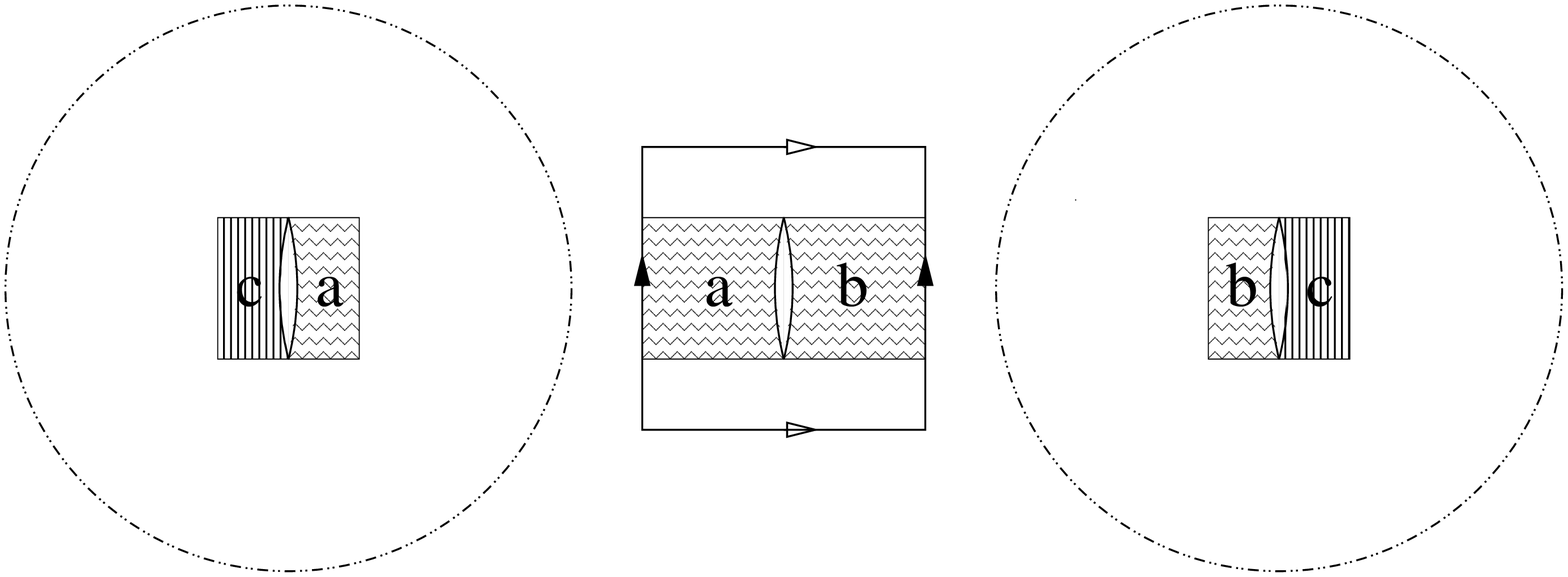}
\caption{Horizontal strips}
\label{fig:strips}
\end{center}
\end{figure}
The width of each strip is $l$, so the sum of the widths of the horizontal strips 
is $2l$. It is interesting that, although the boundary of $\bar Y$ consists of 
three arcs of length $l$, the total vertical length of the boundary of $\bar Y$
is $4l$. This is because $\bar Y$ has two boundary components
of length $2 l$ each, the arc $c$ is traveled twice. Note that the strip 
intersecting arc $c$ would be degenerate if we 
forced the boundaries of strips to land on the boundaries of $\bar Y$.
We set boundaries of strips to be in the expanding annuli
to avoid this problem.
\end{example}

If $\alpha$ is a curve in $X$, a component $Y$ of $X \setminus \alpha$ is 
a subsurface of $X$ whose boundary consists of one or two copies 
of $\alpha$. We consider the annulus whose core curve is $\alpha$ to be 
a component of $X \setminus \alpha$ also. Let $\sigma$ be the hyperbolic 
metric on $X$. 

\begin{theorem} \label{if}
Let $\alpha$ be a curve in $X$, $\beta$ be a transverse curve to
$\alpha$, and $Y$ be a component of $X \setminus \alpha$. 
Assume $l_{\sigma}(\beta) \leq L$. We have:
\begin{enumerate}
\item If $\alpha$ is mostly vertical, then
$  i_Y(\beta, \nup) \prec D_L.$
\item If $\alpha$ is mostly horizontal, then
$  i_Y(\beta, \nm) \prec D_L.$
\end{enumerate}
Here $D_L$ is as in \thmref{q-collar}, with $\log D_L\asymp e^L$.
\end{theorem}

\begin{proof}
Assume $\alpha$ is mostly vertical, that is, 
\begin{equation} \label{vertical}
l_q(\alpha) \leq 2 v_q(\alpha).
\end{equation}
We prove that $i_Y(\beta, \nup) \leq D_L$. From \thmref{q-collar} we have
\begin{equation}
\label{eq1}
l_q(\beta) \leq D_L \, l_q(\alpha).
\end{equation}
Let $\bar \beta$ be the lift of $\beta$ to $\bar X$, the $Y$--cover of $X$. The geometric intersection number of $\bar \beta$ with a horizontal
leaf $h$ equals the intersection number of restrictions of $h$ and 
$\bar \beta$ to $\bar Y$ up to an additive error (\lemref{essential}). 
Let $\bar \beta_0$ be the restriction of
$\bar \beta$ to $\bar Y$. Since $p|_{\bar Y}$ is one-to-one, we have
\begin{equation}
\label{eq2}
l_{\bar q}(\bar \beta_0) \leq l_q (\beta)
\end{equation}
and
\begin{equation} \label{eq3}
v_q(\alpha) \leq v_{\bar q}(\p \bar Y).
\end{equation}
Let $\{ R_1 , \dots , R_p \}$ be the set of maximal horizontal strips. 
Now \eqref{vertical},  \eqref{eq1},  \eqref{eq2}, \eqref{eq3} and 
\corref{strip} imply
\begin{equation}
l_{\bar q}(\bar \beta_0) \leq 4 \, D_L \sum_{i=1}^p w(R_i) .
\end{equation}
Let $m_i$ be the number of times that $\bar \beta_0$ crosses $R_i$ 
and $m$ be the minimum of the $m_i$. Then there exists a horizontal leaf 
$h$ that intersects $\bar \beta_0$ at most $m$ times. That is,
$i_Y(\beta, \nup) \leq m$. But
$$ m \sum_{i=1}^p w(R_i) \leq \sum_{i=1}^p m_i \, w(R_i) 
\leq l_{\bar q}(\bar \beta_0) . $$
Therefore, $m \prec  D_L$. If $\alpha$ is mostly horizontal,
we can similarly show that $i_Y(\beta, \nm) \prec D_L$.
\end{proof}

\begin{theorem} \label{only-if}
Let $\alpha$ be a curve in $X$ and $\beta$ be a transverse curve to
$\alpha$. Assume $l_\sigma(\alpha) \leq \ep$. Then there exists a 
component $Y$ of $X \setminus \alpha$ such that either
\begin{enumerate}
\item $Y$ is an annulus and
$$ \max \big( i_Y(\beta, \nup), i_Y(\beta, \nm) \big) \succ \frac 1\ep, 
\quad {\rm or} $$
\item $Y$ is not an annulus and
$$  \log \Big( \max \big( i_Y(\beta, \nup) \big), i_Y(\beta, \nm) \big)\Big)
\succ \frac 1\ep .$$
\end{enumerate}
\end{theorem}

\begin{proof}
Since the hyperbolic length of $\alpha$ is less than $\ep$, 
\thmref{prim-ann} implies that there exists a primitive annulus
$A$ whose core is in the homotopy class of $\alpha$ and 
whose modulus is $m \asymp \frac 1\ep$. Let $d$ be the 
distance between the boundary components of $A$, and let
$\gamma_0$ be the inner boundary of $A$.  If $A$ is flat, let 
$Y$ be $A$ and k=m; if $A$ is expanding, let $Y$ 
be the component of $X \setminus \alpha$ that contains $A$ 
and let
\begin{equation}  
k= \frac d{ l_q(\gamma_0)}.
\end{equation}
\lemref{p-dist} implies that $\log k \asymp m$.

Let $\bar A$, $\bar \alpha$ and $\bar \beta$ be the lifts of $A$, 
$\alpha_q$ and $\beta_q$ to $\bar X$, respectively, and let $\bar \beta_0$ 
be the restriction of $\bar \beta$ to $\bar Y$. The arc $\bar \beta_0$ has 
to cross $\bar A$; therefore, $l_{\bar q}(\bar \beta_0) \geq d$. Also, 
$\bar \alpha$ is a geodesic; therefore, 
$l_{\bar q}(\bar \alpha) \leq l_q(\gamma_0)$. We have
\begin{equation} \label{bt-k}
\frac{l_{\bar q}(\bar \beta_0)}{l_{\bar q}(\bar \alpha)} \succ k.
\end{equation}
The arc $\bar \beta_0$ is either mostly vertical or mostly horizontal. That is,
\begin{equation} \label{mv}
 2 v_{\bar q}(\bar \beta_0) \geq l_{\bar q}(\bar \beta_0)
\qquad {\rm or} \qquad 
 2 h_{\bar q}(\bar \beta_0) \geq l_{\bar q}(\bar \beta_0).
\end{equation}
Assume the first inequality holds. Let  $\{ R_1, \ldots, R_p \}$ be the 
set of horizontal  strips. \corref{strip}, \eqref{bt-k} and the pigeonhole 
principle imply that $\bar \beta_0$ has to intersect some $R_i$ at least 
(up to a multiplicative constant) $k$ times. But the geometric intersection
number between $\bar \beta$ and the horizontal foliation is equal to
the intersection number of $\bar \beta_0$ with a horizontal leaf up
to an additive constant (\lemref{essential}). That is,
$$  i_Y(\beta, \nup) \succ k .$$
In case the second inequality in \eqref{mv} holds, we consider vertical 
strips and we show similarly that $\bar \beta_0$ intersects a vertical leaf at 
least (up to a multiplicative constant) $k$ times. That is
$$  i_Y(\beta, \nm) \succ k .$$
This finishes the proof of the theorem.
\end{proof}


\section{Proof of the main theorem} \label{sec:proof}

\noindent
The main theorem is a corollary of Theorems \ref{if} and 
\ref{only-if}. First we recall a few facts about the \Teich space. 
The \Teich space of $S$ is the space of conformal structures
on $S$, where two structures are considered to be equivalent if there
is a conformal map between them isotopic to the identity. There are 
several natural metrics defined on $\T(S)$, all inducing the same natural
topology. We work with the Teichm\"uller metric, which assigns to 
$X_1, X_2 \in \T(S)$ the distance
$$d_{\T(S)}(X_1, X_2) = \frac{1}{2} \log(k), $$
where $k$ is the smallest dilatation of a quasi-conformal homeomorphism
from $X_1$ to $X_2$ that is isotopic to the identity. 

Geodesics in the Teichm\"uller space of $S$ are determined by
the quadratic differentials. Let $q$ be a quadratic differential
on a Riemann surface $X \in \T(S)$. Define $g(0)= X$ and, for $t \in \R$,
$g(t)$ to be the conformal structure obtained by scaling the horizontal 
foliation of $q$ by a factor of $e^t$, and the vertical by a factor of 
$e^{-t}$. Then $g\co \R \to \T(S)$ is a geodesic in $\T(S)$ parametrized by 
arc length. The corresponding family of quadratic differentials 
is denoted by $q_t$. For a curve $\alpha$ in $S$, the horizontal and vertical 
lengths of $\alpha$ vary with time as follows: 
\begin{equation}
h_{q_t}(\alpha) = h_q(\alpha) \, e^{-t} \qquad {\rm and} \qquad
v_{q_t}(\alpha) = v_q(\alpha) \, e^t.
\end{equation}
We define a constant $K_{\alpha}$, which measures the relative
complexity of $\mup$ and $\mm$ from the point of view of $\alpha$,
as follows: For each subsurface $Y$ of $X$,
\begin{enumerate}
\item If $Y$ is not an annulus, define $K_Y= log K$ where $K$ is 
the smallest positive number such that $\mup_Y$ and $\mm_Y$ 
have a common $K$--quasi-parallel. 
\item If $Y$ is an annulus, define $K_Y= i_Y(\mup,\mm)$. 
\end{enumerate}
Now, define $K_{\alpha}$ to be the largest $K_Y$ where $Y$ 
is a component of $X \setminus \alpha$.

\begin{theorem} \label{main-long}
Let $(g,\mup,\mm)$ be a geodesic in $\T(S)$ and $\alpha$ be a simple 
closed curve in $S$. Assume $\alpha$ is balanced at $t_0$. Then
$$\frac 1{l_g(\alpha)} \asymp \frac 1{l_{g(t_0)}(\alpha)} \asymp K_\alpha.$$
\end{theorem}

\proof
By definition of $l_g(\alpha)$ we have
\begin{equation}
\frac 1{l_g(\alpha)} \succ \frac 1{l_{g(t_0)}(\alpha)}.
\end{equation}
\thmref{only-if} implies that for every time $t$ there exists
a component $Y$ of $X \setminus \alpha$ such that either
\begin{enumerate}
\item $Y$ is an annulus and $i(\mu^+,\mu^-) \succ 1/ l_{g(t)}(\alpha)$,\footnote{One 
can always find a curve $\beta$ such that
$ i(\mu^+,\mu^-) \geq  \max(i_Y(\beta , \mu^+) , i_Y(\beta, \mu^-))$.} or
\item$Y$ is not an annulus and  $\log(K_Y) \succ 1/l_{g(t)}(\alpha)$.
\end{enumerate}
Therefore, 
\begin{equation}
K_{\alpha} \succ \frac 1{l_g(\alpha)}.
\end{equation}
Let $\ep= l_{g(t_0)}(\alpha)$. It remains to show that, for every component 
$Y$ of $X \setminus \alpha$, $\frac 1\ep \succ K_Y$.

If $Y$ is not an annulus, pick a transverse curve $\beta$ of length less 
than $L$, where $ e^L \asymp \frac 1\ep$. Since $\alpha$ is 
mostly horizontal and mostly vertical simultaneously, \thmref{only-if}
implies that $\beta_Y$ (see \secref{sec:curve} for definition) is 
$D_L$--quasi-parallel to $\mup_Y$ and
$\mm_Y$. Therefore, $K_Y \prec e^L \asymp \frac 1\ep$. 

If $Y$ is an annulus, then $\bar Y$ is a flat annulus. Let $d$ be the
distance between the boundary components of $\bar Y$ and $l$ be
the length of one of the boundary components of $\bar Y$. We know that
$$\frac dl = \Mod(\bar Y) \prec \frac 1\ep.$$
Visualizing $\bar Y$ standing vertically, the restrictions of $\mup$ and 
$\mm$ to $\bar Y$ have angles $\frac \pi 2$ and $- \frac \pi 2$ respectively
(because $\alpha$ is balanced). 
Therefore, an arc in $\mup_Y$ intersects an arc in $\mm_Y$ 
(up to an additive error) $\frac dl$ times. This implies (\lemref{essential})
$$K_Y = i_Y(\mup,\mm) \prec \frac 1\ep.\eqno{\qed}$$

\begin{remark} \label{h=0}
The above theorem is silent with respect to the case when $\alpha$ is 
never balanced (i.e., if the horizontal or the vertical length of $\alpha$ 
is zero along $g$). It is easy to see that, if $\alpha$ is homotopic to a closed 
leaf of $\mup$ or $\mm$, then $l_g(\alpha)=0$, and if not, then $l_g(\alpha)$ 
is positive (consider, for example, the geodesic corresponding to the 
quadratic differential in \exref{example}). However, one cannot approximate 
the value of $l_g(\alpha)$ using the underlying foliations only. 
\end{remark}


 \section{Subsurface distances} \label{sec:cc}

\noindent 
In this section we recall some of the properties of the subsurface 
distances from \cite{CCI} and \cite{CCII} and prove Theorems 
\ref{one}, \ref{two} and \ref{counter}.

First we define a metric on $\CS(S)$ as follows.
For curve systems $\Gamma$ and $\Gamma'$ in $\CS(S)$, define 
$d_S(\Gamma,\Gamma')$ to be equal to one if $[\Gamma] \not = [\Gamma']$ 
and if an element of $[\Gamma]$ and an element of $[\Gamma']$ have 
representatives that are disjoint  from each other. Let the metric on $\CS(S)$ 
be the maximal metric having the above property. That is,
$d_S(\Gamma, \Gamma')= n$ if 
$\Gamma=\Gamma_0, \Gamma_1, \ldots, \Gamma_n=\Gamma'$ is the 
shortest sequence of curve systems on $S$ such that $\Gamma_{i-1}$ is 
distance one from $\Gamma_i$, $i=1, \ldots, n$. 

\begin{theorem}[Masur, Minsky]\label{int-bd}
Let $\Gamma$ and $\Gamma'$ be curve systems in $S$. Then  
$$d_S(\Gamma, \Gamma') \prec \log \big( i_S(\Gamma, \Gamma') \big). $$
\end{theorem}

\begin{remark}
The theorem is true for any subsurface of $S$ that is not an annulus.
If $Y$ is an annulus then
$$d_Y(\Gamma,\Gamma') = i_Y(\Gamma,\Gamma')+1.$$
\end{remark}
Conversely, a bound on all subsurface distances
gives a bound on the intersection number between two curves.

\begin{theorem}[Masur, Minsky]\label{int-sub}
For every $D > 0 $, there exists $K > 0$ such that,
for $\Gamma, \Gamma' \in \CC(S)$, if $i_S(\Gamma,\Gamma') > K$, then 
there exists a subsurface $Y$ of $S$ such that
$$ d_Y(\Gamma, \Gamma') > D.$$
\end{theorem}
\noindent
See \cite{CCI} and \cite{CCII} for the proofs of the above theorems and 
detailed discussion.

Let $\mup$ and $\mm$ be foliations (projectivized 
measured foliations or unmeasured foliations) on $S$. Define 
$d_Y(\mup, \mm)$ to be the distance in $\CS(Y)$ between the restrictions
of $\mup$ and $\mm$ to $Y$ (see \secref{sec:curve} for definition).
We call this the $Y$--distance between $\mup$ and $\mm$.

As we mentioned in the introduction, Minsky has shown in \cite{elcI} that 
every short curve of a hyperbolic 3--manifold $(N, \mup, \mm)$
homeomorphic to $S \times \R$ appears as a  boundary component of 
a subsurface Y of $S$ where the $Y$--distance of $\mup$ and $\mm$ is large.
Now we are ready to prove Theorems \ref{one} and \ref{two}.
For simplicity, we present proofs without keeping track of constants.

\startproof{\thmref{one}}
The hyperbolic 3--manifold $(N,\mup,\mm)$ has bounded geometry 
if and only if there a uniform bound on $d_Y(\mup, \mm)$,
for every subsurface $Y$ of $S$. Also, our main theorem (\thmref{main-long}) 
implies that a geodesic $(g,\mup,\mm)$ has bounded geometry if there 
exists a uniform bound on $K_{\alpha}$, for every curve $\alpha$ in $S$.
We have to show these two conditions are equivalent.

Let $Z$ be a subsurface, $\alpha$ be a boundary component of $Z$ and
$Y$ be a component of $S \setminus \alpha$ containing $Z$. A bound on
$K_\alpha$ provides a bound on $i_Y(\mu^+ , \beta)$ and $i_Y(\beta,
\mu^-)$, where $\beta$ is an arc in $Y$. This provides a bound for
$i_Z(\mu^+ , \beta)$ and $i_Z(\beta, \mu^-)$, which in turn bound
$d_Z(\mu^+, \mu^-)$ (Theorem 7.1). Therefore a uniform bound for
$K_\alpha$ provides a uniform bound for $d_Z(\mu^+, \mu^-)$.

On the other hand, a bound on all subsurface distances gives 
a bound for every $i_Y(\mup, \mm)$ (\thmref{int-sub}). Therefore, 
the two conditions are equivalent. \qed

\startproof{\thmref{two}}
If a curve $\alpha$ is short in $(N, \mu^+, \mu^-)$, then it is a
boundary of a subsurface $Z$, where $d_Z(\mu^+, \mu^-)$ is large. As
we saw above, $K_\alpha(\mu^+, \mu^-)$ has to be large (small
$K_\alpha$ implies small $d_Z(\mu^+, \mu^-)$). Therefore,
$l_g(\alpha)$ has to be small as well.
\endproof

Here we provide a family of examples satisfying \thmref{counter}.

\medskip{\bf \thmref{counter}}\qua {\sl
There exist a sequence of manifolds $(N_n,\mup_n,\mm_n)$,
corresponding Teichm\"uller geodesics $(g_n, \mup_n, \mm_n)$ and 
a curve $\gamma$ in $S$ such that $l_{N_n}(\gamma) \geq c>0$
for all $n$, but
$$l_{g_n}(\gamma) \to 0  \qquad \text{as} \qquad  n \to \infty.$$}

\begin{proof}
Let $(g, \mup, \mm)$ be a geodesic in $\T(S)$, and let 
$\phi_1$ (respectively,  $\phi_2$) be a mapping class of $S$ 
whose support is a subsurface $Y_1$ ($Y_2$) and whose
restriction to $Y_1$ ($Y_2$) is pseudo-Anosov. Assume that
$Y_1$ and $Y_2$ intersect each other nontrivially (that is,
the boundary components of one have nontrivial projections to the 
other), and that they fill a subsurface $Z$ which is a component of
$S \setminus \alpha$, for some curve $\alpha$ (in fact, it would be enough 
to assume $Z$ is a proper subsurface of $S$). Define
$$\mup_n = \phi_1^n(\mup) \qquad  
{\rm and} \qquad \mm_n = \phi_2^n(\mm).$$
Consider the sequence of geodesics $(g_n, \mup_n, \mm_n)$ and
the corresponding sequence of hyperbolic 3--manifolds 
$(N_n, \mup_n, \mm_n)$. Let $K_n$ be the smallest number such that
the restrictions of $\mup_n$ and $\mm_n$ to $Z$ have a common 
$K_n$--quasi-parallel. To prove the theorem, we show that the
length of $\alpha$ remains bounded in $N_n$, but goes to 
zero in $g_n$, as $n$ approaches infinity. It is enough to show that,
as $n$ approaches infinity,
\begin{enumerate}
\item $d_{\bar Z}(\mup_n, \mm_n)$ stays bounded, for any subsurface 
$\bar Z$ which includes $\alpha$ as a boundary component, and
\item $K_n$ approaches infinity.
\end{enumerate}

The curve $\alpha$ is outside of $Y_1$ and $Y_2$. Therefore,
each of $Y_1$ and $Y_2$ either is disjoint from $\bar Z$ or
has a boundary component with has a nontrivial projection to $\bar Z$ 
($ \bar Z \subset Y_1$ is excluded). Assume the restrictions of 
both $\p Y_1$ and $\p Y_2$ to $\bar Z$ are nontrivial.
Since the supports of $\phi_1$ and $\phi_2$ are $Y_1$ and $Y_2$, 
respectively, we have
$$i_{\bar Z}(\mup_n, \p Y_1)   = i_{Z}(\mup, \p Y_1) \qquad  
{\rm and} \qquad i_{\bar Z}(\mm_n, \p Y_2)   = i_{\bar Z}(\mm, \p Y_2).$$
Therefore,
\begin{align*}
d_{\bar Z}(\mup_n, \mm_n)
& \leq d_{\bar Z}(\mup_n, \p Y_1) + d_{\bar Z}(\p Y_1, \p Y_2) 
       + d_{\bar Z}(\p Y_2, \mm_n) \\
& \prec \log i_{\bar Z}(\mup_n, \p Y_1) + d_{\bar Z}(\p Y_1, \p Y_2) 
          + \log i_{\bar Z}(\p Y_2, \mm_n) \\
& =  \log i_{\bar Z}(\mup, \p Y_1) + d_{\bar Z}(\p Y_1, \p Y_2) 
     + \log i_{\bar Z}(\p Y_2, \mm),
\end{align*}
which is independent of $n$. If only one of them, say $Y_1$, is disjoint 
from $\bar Z$, then
$$i_{\bar Z}(\mup_n, \p Y_2)   = i_{\bar Z}(\mup, \p Y_2).$$
Now we can write
$$d_{\bar Z}(\mup_n, \mm_n) \leq 
d_{\bar Z}(\mup_n, \p Y_2) +  d_{\bar Z}(\p Y_2, \mm_n) $$
and argue as before. And finally, if both $Y_1$ and $Y_2$ are 
disjoint from $\bar Z$, then
$$i_{\bar Z}(\mup_n, \mm_n) = i_{\bar Z}(\mup, \mm).$$
This finishes the proof of part (1).

Let $\omega_n$ be the arc that is $K_n$--quasi-parallel to
$\mup_n$ and $\mm_n$. Since $Y_1$ and $Y_2$ fill $Z$,
$\omega_n$ has  to intersect one of them nontrivially. By taking a
subsequence, we can assume that $\omega_n$ intersects $Y_1$
nontrivially (the proof for the subsequence intersecting $Y_2$
nontrivially is similar).
We have
\begin{align*}
d_{Y_1}(\mup_n, \mm_n) 
& \leq d_{Y_1}(\mup_n, \omega_n) +d_{Y_1}(\omega_n, \mm_n) \\
& \prec 2 \log K_n.
\end{align*}
Since $\phi_1$ restricted to $Y_1$ is pseudo-Anosov, we have (see \cite{CCI})
$$d_{Y_1}(\mup_n, \p Y_2) \to \infty \quad {\rm as} \quad n \to \infty.$$
But\vspace{-15pt}
\begin{align*}
 d_{Y_1}(\mm_n, \p Y_2) & \prec \log i_{Y_1}(\mm_n, \p Y_2) \\
                                        & \leq \log i_Z(\mm_n, \p Y_2)\\
                                        & = \log i_Z(\mm, \p Y_2) , \\
\end{align*}\vspace{-30pt}

which is bounded. Therefore, 
$$d_{Y_1}(\mup_n, \mm_n) \to \infty \quad {\rm as} \quad n \to \infty.$$
This implies that $K_n$ approaches infinity as $n$ goes to infinity,
which proves part (2).
\end{proof}



\begin{thebibliography}

\bibitem {bers} {\bf L Bers}, \textit{Simultaneous uniformization}, Bull. Amer.
math. Soc. {66} (1960) 94--97
  \MR{0111834}

\bibitem {bon} {\bf F Bonahon}, \textit{Bouts des vari\'et\'es 
hyperboliques de dimension 3}, Ann. Math. {124} (1986) 71--158
  \MR{0847953}

\bibitem{elcII} {\bf J Brock}, {\bf D Canary}, {\bf Y Minsky},
               \emph{The classification of Kleinian surface groups II:
               the ending lamination conjecture}, in preparation

\bibitem {gard} {\bf F\,P Gardiner}, \textit{Quasiconformal
Teichm\"uller theory}, Mathematical Surveys and Monographs 76,
Amer. Math. Soc. (2000) \MR{1730906}

\bibitem {g-m} {\bf F\,P Gardiner}, {\bf H Masur}, \textit{Extremal
length geometry of Teichm\"uller space}, Complex Variables Theory Appl. {16}
(2000) 209--237 \MR{1099913}

\bibitem {maskit} {\bf B Maskit}, \textit{Comparison of hyperbolic and
extremal lengths}, Ann. Acad. Sci. Fenn. Ser. A I Math. {10}
(1985) 381--386 \MR{0802500}

\bibitem {CCI} {\bf H Masur}, {\bf Y Minsky}, \textit{Geometry of the
complex of curves I, hyperbolicity}, Invent. Math. {138} (1999)
103--149 \MR{1714338}

\bibitem {CCII} {\bf H Masur}, {\bf Y Minsky}, \textit{Geometry of the
complex of curves II, hierarchical structure}, Geom. and
Func. Anal. {10} (2000) 902--974 \MR{1791145}

\bibitem {m-thesis}{\bf Y Minsky}, \textit{Harmonic maps, length, and energy in
Teichm\"uller space}, J. Diff. Geo.  {35} (1992) 151--217
  \MR{1152229}

\bibitem {m-teich} {\bf Y Minsky}, \textit{Teichm\"uller geodesics and ends
of hyperbolic 3--manifolds}, Topology {32} (1993) 625--647
  \MR{1231968}

\bibitem {m-torus} {\bf Y Minsky} \textit{The classification of punctured-torus
group}, Annals of Math. {149} (1999) 559--626
  \MR{1689341}

\bibitem {m-klein} {\bf Y Minsky} \textit{Kleinian groups and the complex of 
curves}, \gtref4{2000}{3}{117}{148}

\bibitem{elcI} {\bf Y Minsky} \emph{The classification of Kleinian
               surface groups I: models and bounds},
               \arxiv{math.GT/0302208}

\bibitem {m-bounded} {\bf Y Minsky} \textit{Bounded geometry for
Kleinian groups}, Invent. Math. {146} (2001) 143--192 \MR{1859020}

\bibitem{thesis} {\bf K Rafi}, \emph{Hyperbolic 3--manifolds and
               geodesics in Teichm\"uller space},  PhD thesis, SUNY
               at Stony Brook (August 2001)
               
\bibitem{rees} {\bf M Rees}, \emph{The geometric model and large
               Lipschitz equivalence direct from Teichm\"uller
               geodesic}, preprint,\nl
               \url{http://www.liv.ac.uk/~maryrees/papershomepage.html}

\bibitem {strebel} {\bf K Strebel}, \textit{Quadratic differentials},
Ergebnisse series 3, Springer--Verlag (1984)
\MR{0743423}

\bibitem {th-note} {\bf W Thurston}, \textit{The Geometry and Topology
of 3--Manifolds}, Princeton University Lecture Notes (1982)
\url{http://www.msri.org/publications/books/gt3m/}

\end{thebibliography}
\end{document}